\documentclass{amsart}
\usepackage{amsfonts}
\usepackage{blkarray}
\usepackage{hyperref,tikz,verbatim,tkz-berge,tkz-graph}
\usepackage{float,graphicx}
\tikzstyle{vertex}=[circle, draw, inner sep=0pt, minimum size=6pt]
\newcommand{\vertex}{\node[vertex]}
\newtheorem{theorem}{Theorem}[section]

\theoremstyle{definition}
\newtheorem{definition}[theorem]{Definition}
\newtheorem{example}[theorem]{Example}
\theoremstyle{remark}
\newtheorem{remark}[theorem]{Remark}
\numberwithin{equation}{section}

\begin{document}
\def\emline#1#2#3#4#5#6{\put(#1,#2){\special{em:moveto}}\put(#4,#5){\special{em:lineto}}}
\def\newpic#1{}

\title[New extremal binary self-dual codes from bicubic planar graphs]{New
extremal binary self-dual codes of lengths 64 and 66 from bicubic planar
graphs}
\author{Ab\.id\.in Kaya}
\address{Department of Computer Engineering, Bursa Orhangazi University,
16310, Bursa, Turkey}
\email{abidin.kaya@bou.edu.tr}
\subjclass[2010]{Primary 94B05, 94B60, 94B65}
\keywords{extremal codes, codes over rings, lifts, bicubic planar graphs}

\begin{abstract}
In this work, connected cubic planar bipartite graphs and related binary
self-dual codes are studied. Binary self-dual codes of length 16 are
obtained by face-vertex incidence matrices of these graphs. By considering
their lifts to the ring $R_2$ new extremal binary self-dual codes of lengths
64 are constructed as Gray images. More precisely, we construct 15 new codes
of length 64. Moreover, 10 new codes of length 66 were obtained by applying
a building-up construction to the binary codes. Codes with these weight
enumerators are constructed for the first time in the literature. The
results are tabulated.
\end{abstract}

\maketitle

\section{Introduction}

Self-dual codes are connected to areas such as design theory, graph theory
and lattices. Therefore they form an interesting family of linear codes.
Such codes have been studied extensively especially the ones over the binary
alphabet. $\left[ 24,12,8\right] _{2}$ extended binary Golay code and $\left[
48,24,12\right] _{2}$ extended quadratic residue code are celebrated
examples of this type. An upper bound for the minimum distance of a binary
self-dual code was given in \cite{conway}. That was finalized in \cite{rains}
as follows; the minimum distance $d$ of a binary self-dual code of length $n$
satisfies $d\leq 4\left[ n/24\right] +6$ if $n\equiv 22\pmod{24}$ and $d\leq
4\left[ n/24\right] +4$, otherwise. A self-dual code is said to be $\emph{%
extremal}$ if it meets the bound. Since the appearance of \cite%
{conway,dougherty,rains} construction and classification of extremal binary
self-dual codes have been a captivating research area.

The existence of extremal binary self-dual codes are open problems for
various lengths. The most famous of these is the existence of doubly-even $%
\left[ 72,36,16\right] _{2}$ self-dual code. The possible weight enumerators
of extremal self-dual binary codes of lengths up to $64$ and $72$ were
determined in \cite{conway}. The weight enumerators for the remaining
lengths up to $100$ were characterized in \cite{dougherty}. Different
techniques such as circulant constructions, Hadamard matrices, automorphism
groups and extensions are used to obtain new extremal binary self-dual
codes. \cite{huffman} is a survey on self-dual codes over different
alphabets. We refer the reader to \cite{dontcheva,dougherty,kayayildiz,kim}
for more information in this direction.

In recent years, new binary self-dual codes have been constructed as Gray
images of self-dual codes over some rings of characteristic $2$. Four
circulant construction was applied to the ring $\mathbb{F}_{2}+u\mathbb{F}%
_{2}$ in\cite{karadeniz}. $\mathbb{F}_{4}+u\mathbb{F}_{4}$-lifts of
self-dual quaternary codes were considered in \cite{kaya}. The ring $\mathbb{%
F}_{2}+u\mathbb{F}_{2}+v\mathbb{F}_{2}+uv\mathbb{F}_{2}$ were used in \cite%
{aksoy,kayayildiz}. In \cite{karadeniz64}, Karadeniz and Yildiz obtained
self-dual codes as lifts of $\left[ 8,4,4\right] _{2}$ extended binary
Hamming code to the ring $R_{3}$ which is of size $2^{8}$.

On the other hand, codes can be obtained from some special graphs by
considering their incidence or adjacency matrices. See \cite{curtis, oral,
oralthesis}.

The possible weight enumerators of extremal singly-even $\left[ 64,32,12%
\right] _{2}$ codes were determined in \cite{conway} as:
\begin{eqnarray*}
W_{64,1} &=&1+\left( 1312+16\beta \right) y^{12}+\left( 22016-64\beta
\right) y^{14}+\cdots ,:14\leq \beta \leq 284, \\
W_{64,2} &=&1+\left( 1312+16\beta \right) y^{12}+\left( 23040-64\beta
\right) y^{14}+\cdots ,:0\leq \beta \leq 277.
\end{eqnarray*}%
The existence of the codes is unknown for the most of the $\beta $ values.
Recently, codes with $\beta =$25, 59 and 74 in $W_{64,1}$ are constructed in
\cite{kayayildiz} by a bordered four circulant construction. Ten new codes
with weight enumerators in $W_{64,2}$ were obtained in \cite{karadeniz64}.
Together with these, codes exist with weight enumerators for $\beta =$14,
18, 22, 25, 29, 32, 36, 39, 44, 46, 53, 59, 60, 64 and 74 in $W_{64,1}$ and
for $\beta =$0, 1, 2, 4, 5,\ 6, 8, 9, 10, 12, 13,\ 14, 16,\ 17, 18, 20, 21,\
22, 23, 24,\ 25,\ 28,\ 29,\ 30, 32,\ 33,\ 36, 37, 38, 40,\ 41,\ 44, 48, 51,\
52,\ 56, 58, 64, 72, 80,\ 88,\ 96, 104, 108,\ 112,\ 114,\ 118,\ 120 and 184
in $W_{64,2}$.

In this work, we obtain fifteen new extremal binary self-dual codes of
length 64. The codes with these weight enumerators are constructed for the
first time in the literature. More precisely, the codes with weight
enumerators for $\beta =$16, 20, 24, 26, 28, 30, 34 and 38 in $W_{64,1};$ $%
\beta =$3, 7, 11, 15, 26, 27 and 35 in $W_{64,2}$ are discovered.

The rest of the work is organized as follows; Section 2 is devoted to the
preliminaries on codes and graphs, bicubic planar graphs and related binary
codes have been studied in Section 3. In Section 4, we consider $R_{2}$%
-lifts of the codes generated by the face-vertex incidence matrices of
connected bicubic planar graphs. As Gray images of these codes extremal
binary self-dual codes of length 64 are constructed. Section 5 is devoted to
the extensions of the new codes of length 64. We were able to obtain ten new
extremal binary self-dual codes of length 66. MAGMA computational algebra
system \cite{magma} have been used for computational results. Section 6
concludes the paper with potential lines of research.

\section{Preliminaries\label{preliminaries}}

The ring $R_{2}=\mathbb{F}_{2}+u\mathbb{F}_{2}+v\mathbb{F}_{2}+uv\mathbb{F}%
_{2}$ was introduced in \cite{yildizr2} as a generalization of $\mathbb{F}%
_{2}+u\mathbb{F}_{2}$. The ring is a commutative local Frobenius ring of
characteristic $2$ and size $16$ that is defined via the restrictions $%
u^{2}=0=v^{2}$ and $uv=vu$. The following isomorphism follows by the
definition
\begin{equation*}
R_{2}\cong \mathbb{F}_{2}[u,v]/\left\langle u^{2},v^{2},uv-vu\right\rangle .
\end{equation*}

Let $R$ denote a commutative Frobenius ring. A linear code $\mathcal{C}$ of
length $n$ over $R$ is an-$R$ submodule of $R^{n}$. The Euclidean dual $%
\mathcal{C}^{\perp }$ of the code $\mathcal{C}$ is defined with respect to
the standard inner product as%
\begin{equation*}
\mathcal{C}^{\perp }:=\left\{ (b_{1},b_{2},\ldots b_{n})\in R^{n}\mid
\sum\limits_{i=1}^{n}a_{i}b_{i}=0,\forall (a_{1},a_{2},\ldots a_{n})\in
\mathcal{C}\right\} .
\end{equation*}%
A code $\mathcal{C}$ is said to be \emph{self-orthogonal} if $\mathcal{C}%
\subseteq \mathcal{C}^{\perp },$ and \emph{self-dual} if $\mathcal{C}=%
\mathcal{C}^{\perp }$. A binary self-dual code is called Type II (\emph{%
doubly-even}) if the weight of any codeword is divisible by 4 and Type I (%
\emph{singly-even}) otherwise.

In \cite{yildizr2} an orthogonality preserving Gray map from $R_{2}^{n}$
into $\mathbb{F}_{2}^{4n}$ was given as follows:
\begin{equation*}
\phi (\bar{a}+u\bar{b}+v\bar{c}+uv\bar{d})=(\bar{d},\bar{c}+\bar{d},\bar{b}+%
\bar{d},\bar{a}+\bar{b}+\bar{c}+\bar{d}),
\end{equation*}%
where $\bar{a},\bar{b},\bar{c},\bar{d}\in \mathbb{F}_{2}^{n}$.

Generator matrices in a special form could be used for lifts. For more
details about the ring $R_{2}$ and lifting a binary code to $R_{2}$ we refer
to \cite{doughertyrk,karadeniz64,yildizr2,aksoy,kayayildiz}.

\begin{definition}
\cite{aksoy} A matrix $\left[ I_{n}|A\right] $ which generates a self-dual
code is called an $LRM$ $\left( \text{lift-ready-matrix}\right) $ if each
upper left $k\times k$ square submatrix of $A$ is invertible.
\end{definition}

Consider the projection $\pi :R_{2}^{n}\rightarrow \mathbb{F}_{2}^{n}$
defined by $\pi \left( \overline{a}+u\overline{b}+v\overline{c}+uv\overline{d%
}\right) =\overline{a},$ where $\overline{a},\overline{b},\overline{c},%
\overline{d}\in \mathbb{F}_{2}^{n}$. Let $\mathcal{D}$ be a binary self-dual
code, a code $\mathcal{C}$ over $R_{2}$ is said to be a \emph{lift} of $%
\mathcal{D}$ if $\pi \left( \mathcal{C}\right) =\mathcal{D}$.

A \emph{graph} $\mathcal{G}=\left( V,E\right) $ is an incidence structure
where $V$ and $E$ are vertices and edges, respectively. A graph that can be
drawn without crossings in the plane is called \emph{planar}. A graph is
said to be \emph{bipartite }if a partition of\ the vertex set $V$ into $%
\left\{ V_{1},V_{2}\right\} $ such that there is no edge with both endpoints
in $V_{1}$ or $V_{2}$.

The graph $\mathcal{G}$ is called \emph{simple} if it is free of loops and
multi-edges. In a simple graph, the degree of a vertex $w$ denoted by $\deg
\left( w\right) $ is the number of edges incident with $w$. A graph is
called $k$-$\emph{regular}$ when every vertex has degree $k$. A $3$-regular
graph is said to be \emph{cubic}. A \emph{bicubic} graph is both bipartite
and cubic. If there is a path between every pair of vertices than it is said
to be \emph{connected}.

\section{Bicubic planar graphs and binary self-dual codes}

In this section, we consider connected bicubic planar graphs and related
self-dual codes. With respect to \cite{oral} any cubic planar graph could be
used to construct self-orthogonal codes but we focus on the following result
for bipartite graphs which gives self-dual codes.

\begin{theorem}
\label{selfdual} \cite{oral} Let $\mathcal{G}$ be a connected cubic planar
bipartite graph with vertex set $\left\{ 1,2,\ldots ,n\right\} $ and
face-vertex incidence matrix $D$. Let $f_{1}$, $f_{2}$ be any two faces of $%
\mathcal{G}$ of different colours in a $3$-face coloring of $\mathcal{G}$.
If we delete the rows corresponding to $f_{1}$ and $f_{2}$ from $D$, the
resulting matrix is a generator matrix for a self-dual code of length $n$.
Moreover, this code is independent of the choice of faces $f_{1}$, $f_{2}$.
\end{theorem}

The extended binary Hamming code of length $8$ is obtained from the
face-vertex incidence matrix of the cube in the following example:

\begin{example}
Consider the cube which is a connected bicubic planar graph.
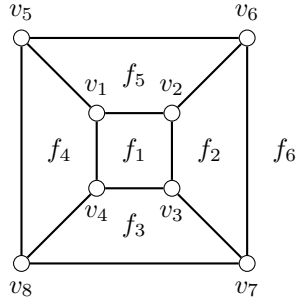
\begin{center}
\GraphInit[vstyle=Hasse]%
\begin{tikzpicture}\centering  [scale=1.2]
        \vertex[label=$v_1$](v1) at (1,2) {};
        \vertex[label=$v_2$](v2) at (2,2) {};
        \vertex[label=below:$v_3$](v3) at (2,1) {};
        \vertex[label=below:$v_4$](v4) at (1,1) {};
        \vertex[label=$v_5$](v5) at (0,3) {};
        \vertex[label=$v_6$](v6) at (3,3) {};
        \vertex[label=below:$v_7$](v7) at (3,0) {};
        \vertex[label=below:$v_8$](v8) at (0,0) {};
    \tikzset{EdgeStyle/.style={-}}
        \Edge(v1)(v2)
        \Edge(v1)(v4)
        \Edge(v1)(v5)
        \Edge(v2)(v6)
        \Edge(v2)(v3)
        \Edge(v3)(v7)
        \Edge(v3)(v4)
        \Edge(v4)(v8)
        \Edge(v5)(v6)
        \Edge(v6)(v7)
        \Edge(v7)(v8)
        \Edge(v5)(v8);
        \draw(1.5,1.5) node{$f_1$};
        \draw(2.5,1.5) node{$f_2$};
        \draw(1.5,0.5) node{$f_3$};
        \draw(0.5,1.5) node{$f_4$};
        \draw(1.5,2.5) node{$f_5$};
        \draw(3.5,1.5) node{$f_6$};
    \end{tikzpicture}
\end{center}
Obviously, the faces $f_{5}$ and $f_{6}$ will have diffterent colors in a $3$%
-face coloring of the cube. Hence, by Theorem \ref{selfdual} the following
submatrix of the face-vertex incidence matrix of the graph generates a
self-dual code.
\begin{equation*}
\begin{array}{c|cccccccc}
& v_{1} & v_{2} & v_{3} & v_{4} & v_{5} & v_{6} & v_{7} & v_{8} \\ \hline
f_{1} & 1 & 1 & 1 & 1 & 0 & 0 & 0 & 0 \\
f_{2} & 0 & 1 & 1 & 0 & 0 & 1 & 1 & 0 \\
f_{3} & 0 & 0 & 1 & 1 & 0 & 0 & 1 & 1 \\
f_{4} & 1 & 0 & 0 & 1 & 1 & 0 & 0 & 1%
\end{array}%
\end{equation*}%
The code is the extended binary Hamming code of parameters $\left[ 8,4,4%
\right] _{2}$. When we express the matrix in standard form, we observe that
it is a lift-ready matrix in the form;%
\begin{equation*}
\left[
\begin{array}{cccc|cccc}
1 & 0 & 0 & 0 & 1 & 0 & 1 & 1 \\
0 & 1 & 0 & 0 & 0 & 1 & 1 & 1 \\
0 & 0 & 1 & 0 & 1 & 1 & 1 & 0 \\
0 & 0 & 0 & 1 & 1 & 1 & 0 & 1%
\end{array}%
\right]
\end{equation*}%
In \cite{karadeniz64}, lifts of the extended binary Hamming code to the ring
$R_{3}$ have been studied in detail, where $R_{3}:=\mathbb{F}%
_{2}[u_{1},u_{2},u_{3}]/\left\langle
u_{i}^{2},u_{i}u_{j}-u_{j}u_{i}\right\rangle $, $1\leq i,j\leq 3$. They were
able to obtain ten new extremal binary self-dual codes with weight
enumerators in $W_{64,2}$.
\end{example}

In the following we give an example that is obtained by combining
two cycles of length $8$. That is the graph E16 in
\cite{oralthesis}.

\begin{figure}[tbph]
\centering
\begin{tikzpicture}[scale=1]
        \vertex(v1) at (1.4,0) {};
        \vertex(v2) at (3.6,0) {};
        \vertex(v3) at (5,1.4) {};
        \vertex(v4) at (5,3.6) {};
        \vertex(v5) at (3.6,5) {};
        \vertex(v6) at (1.4,5) {};
        \vertex(v7) at (0,3.6) {};
        \vertex(v8) at (0,1.4) {};
        \vertex(v9) at (1.9,1) {};
        \vertex(v10) at (3.1,1) {};
        \vertex(v11) at (4,1.9) {};
        \vertex(v12) at (4,3.1) {};
        \vertex(v13) at (3.1,4) {};
        \vertex(v14) at (1.9,4) {};
        \vertex(v15) at (1,3.1) {};
        \vertex(v16) at (1,1.9) {};
    \tikzset{EdgeStyle/.style={-}}
        \Edge(v1)(v2)
        \Edge(v2)(v3)
        \Edge(v3)(v4)
        \Edge(v4)(v5)
        \Edge(v5)(v6)
        \Edge(v6)(v7)
        \Edge(v7)(v8)
        \Edge(v8)(v1)
        \Edge(v9)(v10)
        \Edge(v10)(v11)
        \Edge(v11)(v12)
        \Edge(v12)(v13)
        \Edge(v13)(v14)
        \Edge(v14)(v15)
        \Edge(v15)(v16)
        \Edge(v16)(v9)
        \Edge(v2)(v10)
        \Edge(v3)(v11)
        \Edge(v4)(v12)
        \Edge(v5)(v13)
        \Edge(v6)(v14)
        \Edge(v7)(v15)
        \Edge(v8)(v16)
        \Edge(v1)(v9);
    \end{tikzpicture}
\caption{The graph $\mathcal{G}_{1}$}
\label{fig:g1}
\end{figure}
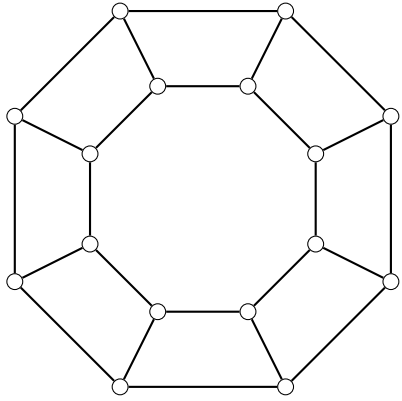

\begin{example}
\label{exg1}Let us consider the graph $\mathcal{G}_{1}$ in Figure \ref%
{fig:g1}. Let $\mathcal{C}$ be the code generated by the face-vertex
incidence matrix of the graph. Then, $\mathcal{C}$ is a self-dual binary
code of length $16$ by Theorem \ref{selfdual}. A generator matrix in
standard form is $\left[ I_{8}|A_{1}\right] ,$ where
\begin{equation}
A_{1}=\left[
\begin{array}{cccccc|cc}
1 & 0 & 0 & 0 & 0 & 0 & 1 & 1 \\
0 & 1 & 0 & 0 & 0 & 0 & 1 & 1 \\
0 & 0 & 1 & 0 & 0 & 0 & 1 & 1 \\
0 & 0 & 0 & 1 & 0 & 0 & 1 & 1 \\
0 & 0 & 0 & 0 & 1 & 0 & 1 & 1 \\
0 & 0 & 0 & 0 & 0 & 1 & 1 & 1 \\ \hline
1 & 1 & 1 & 1 & 1 & 1 & 1 & 0 \\
1 & 1 & 1 & 1 & 1 & 1 & 0 & 1%
\end{array}%
\right]  \label{g1matrix}
\end{equation}%
The code $\mathcal{C}$ is a self-dual Type II code of parameters $\left[
16,8,4\right] _{2}$ with an automorphism group of order $2^{14}3^{2}5\times
7 $ and weight distribution $1+28z^{4}+198z^{8}+28z^{12}+z^{16}$.
\end{example}

\section{Lifts of self-dual graph codes\label{lifts}}

Lifts of the binary self-dual codes related to some bicubic planar graphs on
16 vertices have been considered in this section. $\mathbb{F}_{2}+u\mathbb{F}%
_{2}$-lifts of four circulant binary self-dual codes of length $32$ have
been studied in \cite{karadeniz}. Lifting quaternary codes to $\mathbb{F}%
_{4}+u\mathbb{F}_{4}$ have been used in \cite{kaya}. For more research in
this direction we refer the reader to \cite%
{karadeniz,karadeniz64,aksoy,kaya,kayayildiz,tufekci}. The binary self-dual
codes obtained from face-vertex incidence matrix of connected bicubic planar
graphs is in good shape which is suitable for lifting. In this section, by
considering binary images of $R_{2}$-lifts of graph codes we obtain 15 new
extremal binary self-dual Type I codes of length 64. The existence of
extremal self-dual codes of length 64 with these weight enumerators was
previosly unknown.

We need a brief representation for the elements of the ring $R_{2}$ in order
to tabulate the results. We prefer to use hexadecimals. The correspondence
between the binary $4$ tuples and the hexadecimals is as follows:
\begin{eqnarray*}
0 &\leftrightarrow &0000,\ 1\leftrightarrow 0001,\ 2\leftrightarrow 0010,\
3\leftrightarrow 0011, \\
4 &\leftrightarrow &0100,\ 5\leftrightarrow 0101,\ 6\leftrightarrow 0110,\
7\leftrightarrow 0111, \\
8 &\leftrightarrow &1000,\ 9\leftrightarrow 1001,\ A\leftrightarrow 1010,\
B\leftrightarrow 1011, \\
C &\leftrightarrow &1100,\ D\leftrightarrow 1101,\ E\leftrightarrow 1110,\
F\leftrightarrow 1111.
\end{eqnarray*}%
The ordered basis $\left\{ uv,v,u,1\right\} $ is used to express elements of
$R_{2}$ For instance, $1+u+v$ is represented as $0111$ which corresponds to
hexadecimal $7$.

\begin{example}
The generator matrix $\left[ I_{8}|A_{1}\right] $ of the code
$\mathcal{C}$ in Example \ref{exg1} that is obtained from the graph
$\mathcal{G}_{1}$ is in lift-ready form. $A_{1}$ is lifted to
$K_{1}$ where
\begin{equation*}
K_{1}=\left[
\begin{array}{cccccccc}
9 & C & 0 & 8 & E & 4 & D & 7 \\
4 & 5 & 4 & E & 8 & 8 & B & 1 \\
E & 2 & 1 & 6 & 2 & C & F & B \\
E & 8 & 6 & 9 & 6 & A & F & 7 \\
4 & A & E & A & 1 & A & F & 7 \\
8 & C & 6 & 0 & C & B & 3 & 5 \\
7 & 5 & F & 7 & 3 & B & 5 & 8 \\
B & 9 & D & 5 & D & B & E & 5%
\end{array}%
\right] .
\end{equation*}%
Let $\mathcal{K}_{1}$ be the code over $R_{2}$ generated by $\left[
I_{8}|K_{1}\right] $ then its Gray image $\phi \left(
\mathcal{K}_{1}\right) $ is a self-dual Type I $\left[
64,32,12\right] _{2}$-code with weight enumerator for $\beta =20$ in
$W_{64,1}$. This the first example of an extremal Type I code of
length $64$ with this weight enumerator.
\end{example}

As lifts of $\mathcal{C}$ we were able to obtain five new extremal binary
self-dual Type I codes. We were able to construct codes with weight \
enumerators $\beta =20,24,26$ and $30$ in $W_{64,1}$; $\beta =26$ in $%
W_{64,2}$. The codes are given below in Table \ref{tab:table1}. In
order to save space only the upper triangular parts of the matrices
are given since the rest is determined by orthogonality relations.
\begin{table}[tbph]
\caption{New extremal self-dual Type I $\left[ 64,32,12\right] _{2}$ codes
as binary images of $R_{2}$-lifts of $\mathcal{C}$}
\label{tab:table1}%
\begin{tabular}{||c|c||c||}
\hline $\mathcal{K}_{i}$ & Upper triangular part of the matrix $K_i$
& $\beta $ in $W_{64,i}$
\\ \hline\hline
$\mathcal{K}_{1}$ & $9C08E4D754E88B1162CFB96AF71AF7B35585$ & \textbf{20 }in $%
W_{64,1}$ \\ \hline
$\mathcal{K}_{2}$ & $9A8C663FF2A4855D2463516C7D943F95BB8B$ & \textbf{24 }in $%
W_{64,1}$ \\ \hline
$\mathcal{K}_{3}$ & $D24022373664C1D1AC671120799C7759FB0F$ & \textbf{26 }in $%
W_{64,1}$ \\ \hline
$\mathcal{K}_{4}$ & $9E4CEEF332A0C55D6CE755A83D5C3FD9F78B$ & \textbf{30 }in $%
W_{64,1}$ \\ \hline
$\mathcal{K}_{5}$ & $9A8C6273FEECC119A8E75D6CF51CF3513F87$ & \textbf{26 }in $%
W_{64,2}$ \\ \hline
\end{tabular}%
\end{table}

\begin{remark}
The extremal binary self-dual codes in Table \ref{tab:table1} are
constructed by considering the binary image $\phi \left( \mathcal{K}%
_{i}\right) $ where $\mathcal{K}_{i}$ is the code of length $16$ over $R_{2}$
that is generated by $\left[ I_{8}|K_{i}\right] $.
\end{remark}

Let us recall the graph F16 from \cite{oralthesis} which is a
connected bicubic planar graph. Its face-vertex incidence matrix
leads to a Type I self-dual code of length $16$:
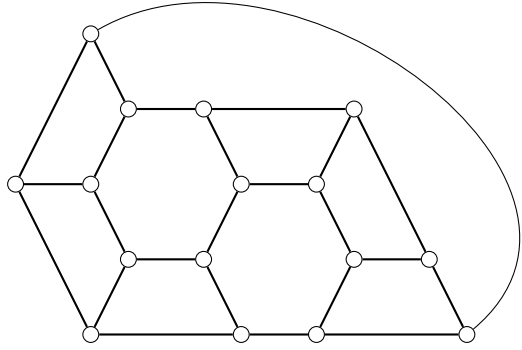
\begin{figure}[tbph]
\centering
\begin{tikzpicture}[scale=1]
        \vertex(v1) at (1.5,1) {};
        \vertex(v2) at (2.5,1) {};
        \vertex(v3) at (3,2) {};
        \vertex(v4) at (2.5,3) {};
        \vertex(v5) at (1.5,3) {};
        \vertex(v6) at (1,2) {};
        \vertex(v7) at (0,2) {};
        \vertex(v8) at (1,4) {};
        \vertex(v9) at (1,0) {};
        \vertex(v10) at (3,0) {};
        \vertex(v11) at (4,0) {};
        \vertex(v12) at (4.5,1) {};
        \vertex(v13) at (4,2) {};
        \vertex(v14) at (4.5,3) {};
        \vertex(v15) at (5.5,1) {};
        \vertex(v16) at (6,0) {};
    \tikzset{EdgeStyle/.style={-}}
        \Edge(v1)(v2)
        \Edge(v2)(v3)
        \Edge(v3)(v4)
        \Edge(v4)(v5)
        \Edge(v5)(v6)
          \Edge(v5)(v8)
        \Edge(v6)(v1)
        \Edge(v7)(v8)
        \Edge(v7)(v9)
        \Edge(v7)(v6)
        \Edge(v1)(v9)
        \Edge(v9)(v10)
        \Edge(v2)(v10)
        \Edge(v10)(v11)
        \Edge(v11)(v12)
        \Edge(v12)(v13)
        \Edge(v13)(v3)
        \Edge(v13)(v14)
         \Edge(v4)(v14)
         \Edge(v11)(v16)
        \Edge(v15)(v16)
        \Edge(v12)(v15)
         \Edge(v14)(v15)
         \draw [black]   (v8) to[out=-330,in=40, distance=3cm
         ](v16);
        ;
    \end{tikzpicture}
\caption{The graph $\mathcal{G}_{2}$}
\label{fig:F16}
\end{figure}
The graph $\mathcal{G}_{2}$ in Figure \ref{fig:F16} is a connected biqubic
planar graph. So, by Theorem \ref{selfdual} the face-vertex incidence matrix
of $\mathcal{G}_{2}$ generates a self-dual binary code $\mathcal{D}$. The
code is Type I with parameters $\left[ 16,8,4\right] _{2}$. The generator
matrix in standard form is $\left[ I_{8}|A_{2}\right] $where:
\begin{equation}
A_{2}=\left[
\begin{array}{ccc|ccc|cc}
1 & 0 & 0 & 1 & 1 & 1 & 1 & 0 \\
0 & 1 & 0 & 1 & 1 & 1 & 1 & 0 \\
0 & 0 & 1 & 1 & 1 & 1 & 1 & 0 \\ \hline
0 & 0 & 0 & 1 & 0 & 0 & 1 & 1 \\
0 & 0 & 0 & 0 & 1 & 0 & 1 & 1 \\
0 & 0 & 0 & 0 & 0 & 1 & 1 & 1 \\ \hline
1 & 1 & 1 & 1 & 1 & 1 & 0 & 1 \\
1 & 1 & 1 & 0 & 0 & 0 & 1 & 1%
\end{array}%
\right]  \label{g2matrix}
\end{equation}%
Weight distribution of $\mathcal{D}$ is $%
1+12z^{4}+64z^{6}+102z^{8}+64z^{10}+12z^{12}+z^{16}$ and $\left\vert
Aut\left( \mathcal{D}\right) \right\vert =2^{13}3^{2}$. The matrix is an $%
LRM $ and as binary images of the $R_{2}$-lifts of $\mathcal{D}$ we obtain
fifteen Type I $\left[ 64,32,12\right] _{2}$-codes with previously unknown
weight enumerators. More precisely, the codes with weight enumerators $\beta
=$16, 20, 24, 26, 28, 30, 34 and 38 in $W_{64,1}$; $\beta =$3, 7, 11, 15,
26, 27 and 35 in $W_{64,2}$ $\ $Those are listed in Table \ref{tab:table2}.

\begin{table}[tbph]
\caption{New extremal self-dual Type I $\left[ 64,32,12\right] _{2}$ codes
as binary images of $R_{2}$-lifts of $\mathcal{D}$ (15 new codes)}
\label{tab:table2}%
\begin{tabular}{||c|c||c||}
\hline $\mathcal{L}_{i}$ & Upper triangular part of the matrix $L_i$
& $\beta $ in $W_{64,i}$
\\ \hline\hline
$\mathcal{L}_{1}$ & $5EA3BBDE945739ADDBB3436ABF7237B5105B$ & \textbf{16 }in $%
W_{64,1}$ \\ \hline
$\mathcal{L}_{2}$ & $F4ED975832719FA15953274813383F97DC7F$ & \textbf{20 }in $%
W_{64,1}$ \\ \hline
$\mathcal{L}_{3}$ & $DAE7379AD85FBDE1D3BB0BAA33F6FFF5D0DF$ & \textbf{24 }in $%
W_{64,1}$ \\ \hline
$\mathcal{L}_{4}$ & $5A2B7796DC53F1E99FBBCFAAB77B37F1D097$ & \textbf{26 }in $%
W_{64,1}$ \\ \hline
$\mathcal{L}_{5}$ & $5E2F335610D7FDA557FF0BAEF3F237F59813$ & \textbf{28 }in $%
W_{64,1}$ \\ \hline
$\mathcal{L}_{6}$ & $D6A3BBDED05739A99BB34FAABF323F3D1C13$ & \textbf{30 }in $%
W_{64,1}$ \\ \hline
$\mathcal{L}_{7}$ & $3C211FD8B23D1FAD115F670C5BB83F9F94F3$ & \textbf{34 }in $%
W_{64,1}$ \\ \hline
$\mathcal{L}_{8}$ & $5205997E1A77550739D9AA92817F01B5358B9$ & \textbf{38 }in
$W_{64,1}$ \\ \hline
$\mathcal{L}_{9}$ & $742D979876F15F2599936F041BFCBF135437$ & \textbf{3 }in $%
W_{64,2}$ \\ \hline
$\mathcal{L}_{10}$ & $16055DB29A3B990771952D6013B09F53D0B5$ & \textbf{7 }in $%
W_{64,2}$ \\ \hline
$\mathcal{L}_{11}$ & $3C295F50FA399B619D9F6F481BF83B57D8BF$ & \textbf{11} in
$W_{64,2}$ \\ \hline
$\mathcal{L}_{12}$ & $56491536127F9147B55DE5241FF0D757587D$ & \textbf{15} in
$W_{64,2}$ \\ \hline
$\mathcal{L}_{13}$ & $5A23B31AD8177DA55B7BC7A6BB7AFF71981B$ & \textbf{26 }in
$W_{64,2}$ \\ \hline
$\mathcal{L}_{14}$ & $DEC1D9F2D63FD94B9D5A12CD330D35B1C3D$ & \textbf{27 }in $%
W_{64,2}$ \\ \hline
$\mathcal{L}_{15}$ & $1023DB7472337DAD9995ED485DB2D3715457$ & \textbf{35 }in
$W_{64,2}$ \\ \hline
\end{tabular}%
\end{table}

\begin{remark}
The binary codes in Table \ref{tab:table2} are constructed as the Gray image
$\phi \left( \mathcal{L}_{i}\right) $ where $\mathcal{L}_{i}$ is the code of
length $16$ over $R_{2}$ that is generated by $\left[ I_{8}|L_{i}\right] $.
\end{remark}

The $\beta $ parameters of the codes in Table \ref{tab:table1} reoccurs in
Table \ref{tab:table2}. We see that the corresponding codes are not
equivalent when we check the invariants.\ If two extremal binary self-dual
Type I codes of length 64 have the same weight enumerator, then for each
code let $c_{1},c_{2},\ldots ,c_{N}$ be the codewords of weight $12$. Let $%
A_{j}=\left\vert \left\{ \left( c_{k},c_{l}\right) |\ d\left(
c_{k},c_{l}\right) =j,\ k<l\right\} \right\vert $ where $d$ is the Hamming
distance. $A_{12}$ is invariant under a permutation of the coordinates.
Hence, two codes are inequivalent if their $A_{12}$-values are not equal.
Those are given in Table \ref{tab:equivalence} which indicates that the
corresponding codes are not equivalent.

\begin{table}[tbph]
\caption{The inequivalence of the codes in Table \protect\ref{tab:table1}
and Table \protect\ref{tab:table2}}
\label{tab:equivalence}%
\begin{tabular}{||c|c||c|c||c||}
\hline
$\mathcal{K}_{i}$ & $A_{12}$ & $\mathcal{L}_{i}$ & $A_{12}$ & $\beta $ in $%
W_{64,i}$ \\ \hline\hline
$\mathcal{K}_{1}$ & $15732$ & $\mathcal{L}_{2}$ & $14964$ & \textbf{20 }in $%
W_{64,1}$ \\ \hline
$\mathcal{K}_{2}$ & 16488 & $\mathcal{L}_{3}$ & $17264$ & \textbf{24 }in $%
W_{64,1} $ \\ \hline
$\mathcal{K}_{3}$ & $17676$ & $\mathcal{L}_{4}$ & $17898$ & \textbf{26 }in $%
W_{64,1}$ \\ \hline
$\mathcal{K}_{4}$ & 20544 & $\mathcal{L}_{6}$ & $19890$ & \textbf{30 }in $%
W_{64,1} $ \\ \hline
$\mathcal{K}_{5}$ & 18876 & $\mathcal{L}_{13}$ & $19680$ & \textbf{26 }in $%
W_{64,2}$ \\ \hline
\end{tabular}%
\end{table}

\begin{theorem}
The existence of extremal binary self-dual codes of length $64$ is known for
$23$ parameters in $W_{64,1};$ $56$ parameters in $W_{64,2}$.
\end{theorem}

\begin{remark}
The codes in Table \ref{tab:table1} and Table \ref{tab:table2} have an
automorphism group of order $2^{2}$. The $R_{2}$ and binary generator
matrices of these are available online at \cite{web}. We also note that the
symmetry in the graphs $\mathcal{G}_{1}$ and $\mathcal{G}_{2}$ emerge in the
standard forms of the generator matrices of the corresponding binary
self-dual codes.
\end{remark}

\section{New binary self-dual codes of length $66$ as extensions\label%
{newcodes}}

Building-up construction which is also known as extension in the literature
is an efficient method to construct self-dual codes from shorther ones. We
refer to \cite{kim,doughertyfrobenius} for different versions of the
construction. Such methods have been effectively used recently in \cite%
{karadeniz,aksoy,kaya,kayayildiz,tufekci} to obtain new self-dual codes of
lenghts 58, 66 and 68. In this section, we apply the following extension
method to the codes in Section \ref{lifts}. As a result of this, 10 new
extremal binary self-dual codes of length 66 are obtained.

\begin{theorem}
$($\cite{doughertyfrobenius}$)$ \label{ext}Let $\mathcal{C}$ be a self-dual
code over $R$ of length $n$ and $G=(r_{i})$ be a $k\times n$ generator
matrix for $\mathcal{C}$, where $r_{i}$ is the $i$-th row of $G$, $1\leq
i\leq k$. Let $c$ be a unit in $R$ such that $c^{2}=1$ and $X$ be a vector
in $R^{n}$ with $\left\langle X,X\right\rangle =1$. Let $y_{i}=\left\langle
r_{i},X\right\rangle $ for $1\leq i\leq k$. Then the following matrix%
\begin{equation*}
\left(
\begin{array}{cc|c}
1 & 0 & X \\ \hline
y_{1} & cy_{1} & r_{1} \\
\vdots & \vdots & \vdots \\
y_{k} & cy_{k} & r_{k}%
\end{array}%
\right) ,
\end{equation*}%
generates a self-dual code $\mathcal{C}^{\prime }$ over $R$ of length $n+2$.
\end{theorem}

A self-dual $\left[ 66,33,12\right] _{2}$-code has a weight enumerator in
one of the following forms (\cite{dougherty})%
\begin{eqnarray*}
W_{66,1} &=&1+\left( 858+8\beta \right) y^{12}+\left( 18678-24\beta \right)
y^{14}+\cdots \text{ where }0\leq \beta \leq 778, \\
W_{66,2} &=&1+1690y^{12}+7990y^{14}+\cdots \text{ } \\
\text{and }W_{66,3} &=&1+\left( 858+8\beta \right) y^{12}+\left(
18166-24\beta \right) y^{14}+\cdots \text{ where }14\leq \beta \leq 756,
\end{eqnarray*}%
Recently, new codes with weight enumerators in $W_{66,1}$ are constructed in
\cite{karadeniz,kaya,tufekci}. More precisely, 5 new codes are obtained in
\cite{karadeniz}, 24 new codes in \cite{kaya} and 11 new codes in \cite%
{tufekci}. Together with these, the existence of such codes is known for $%
\beta =$0, 1, 2, 3, 5, 6, 8$,\ldots ,$11, 14$,\ldots ,$56, 59$,\ldots ,$69,
71$,\ldots ,$ 90, 92, 94 and 100 in $W_{66,1}$.

We construct the codes with weight enumerators $\beta =$13 and 57 in $%
W_{66,1}$.

Most recently, 14 codes were discovered in \cite{kayayildiz} by
applying the building-up construction to the binary images of
modified four-circulant codes of length $16$ over $R_2$. Together
with these the existence of codes in $W_{66,3}$ is known for $\beta
=$28, 29, 30, 31, 32, 33, 34, 35, 36, 37, 38, 43, 44, 45, 46, 47,
48, 49, 50, 51, 52, 53, 54, 55, 56, 57, 58, 59, 60, 61, 62, 63, 64,
64, 66, 67, 70, 71, 73, 74, 75, 76, 77, 78, 79, 80, 81, 82, 83, 84,
85, 86, 87, 88, 90, 92.

In this work, eight new codes with new weight enumerators are constructed.
More precisely the codes with weight enumerator for $\beta =$24, 25, 26, 27,
39, 40, 41 and 42 in $W_{66,3}$ which are listed in Table \ref{tab:table3}.

\begin{center}
\begin{table}[tbph]
\caption{10 new extremal self-dual binary codes of length $66$ by Theorem
\protect\ref{ext}}
\label{tab:table3}%
\begin{tabular}{||c|c|c||c||}
\hline
$\mathcal{C}_{i}$ & Code & $X$ & $\beta $ in $W_{66,i}$ \\ \hline\hline
$\mathcal{C}_{1}$ & $\mathcal{L}_{9}$ &
\begin{tabular}{l}
$00100110100110000010100011011000$ \\
$00100101110100001100000011110001$%
\end{tabular}
& \textbf{13 }in $W_{66,1}$ \\ \hline
$\mathcal{C}_{2}$ & $\mathcal{L}_{15}$ & $01101100110101100100110101100011%
\mathbf{1}^{32}$ & \textbf{57 }in $W_{66,1}$ \\ \hline
$\mathcal{C}_{3}$ & $\mathcal{L}_{1}$ & $00010000111100101101111111100001%
\mathbf{1}^{32}$ & \textbf{24 }in $W_{66,3}$ \\ \hline
$\mathcal{C}_{4}$ & $\mathcal{L}_{1}$ &
\begin{tabular}{l}
$01111110101011011111000101000011$ \\
$00111110001000011000101001110100$%
\end{tabular}
& \textbf{25 }in $W_{66,3}$ \\ \hline
$\mathcal{C}_{5}$ & $\mathcal{K}_{1}$ & $00100011001001110101011010001101%
\mathbf{1}^{32}$ & \textbf{26 }in $W_{66,3}$ \\ \hline
$\mathcal{C}_{6}$ & $\mathcal{L}_{1}$ & $11011000100101000100111110110110%
\mathbf{0}^{32}$ & \textbf{27 }in $W_{66,3}$ \\ \hline
$\mathcal{C}_{7}$ & $\mathcal{K}_{4}$ & $01011011011101111111000111011100%
\mathbf{1}^{32}$ & \textbf{39 }in $W_{66,3}$ \\ \hline
$\mathcal{C}_{8}$ & $\mathcal{L}_{4}$ & $00110110101011010110001010110110%
\mathbf{1}^{32}$ & \textbf{40} in $W_{66,3}$ \\ \hline
$\mathcal{C}_{9}$ & $\mathcal{K}_{3}$ & $00100110110111100100010101111001%
\mathbf{1}^{32}$ & \textbf{41 }in $W_{66,3}$ \\ \hline
$\mathcal{C}_{10}$ & $\mathcal{K}_{3}$ & $10101101100111011001110110010101%
\mathbf{1}^{32}$ & \textbf{42 }in $W_{66,3}$ \\ \hline
\end{tabular}%
\end{table}
\end{center}

\begin{remark}
The extremal codes of length $66$ in Table \ref{tab:table3} are obtained by
considering the binary code generated by the matrix%
\begin{equation*}
\left(
\begin{array}{cc|c}
1 & 0 & X \\ \hline
y_{1} & cy_{1} & r_{1} \\
\vdots & \vdots & \vdots \\
y_{k} & cy_{k} & r_{k}%
\end{array}%
\right)
\end{equation*}%
where $G=(r_{i})$ is the matrix determined by the Gray images of $\left[
I_{8}|L_{i}\right] ,\ u\left[ I_{8}|L_{i}\right] ,\ v\left[ I_{8}|L_{i}%
\right] $ and $uv\left[ I_{8}|L_{i}\right] $. The binary generator matrices
of the codes in Table \ref{tab:table3} are available online at \cite{web}.
The codes all have an automorphism group of order $2$.
\end{remark}

\begin{theorem}
The existence of extremal binary self-dual codes of length $66$ is known for
$89$ parameters in $W_{66,1};$ $64$ parameters in $W_{66,3}$.
\end{theorem}

\section{Conclusion}

The codes obtained by considering face-vertex incidence matrices of bicubic
planar graphs have a nice structure. In this work, we considered two such
graphs on $16$ vertices. New extremal binary self-dual codes of length $64$
obtained as the Gray images of $R_{2}$-lifts of face-vertex incidence
matrices of the graphs. Such methods can be applied to different bicubic
planar graphs and lifts can be considered over various rings. Another
possible research area is to generalize Oral's work to a larger family of
planar graphs.

\end{document}